\documentclass[12pt]{article}
\usepackage{graphicx}
\usepackage{amssymb,amsfonts,amsthm}
\usepackage{amsmath,amsthm,amssymb}
\usepackage{amsfonts}

\newtheorem{theorem}{Theorem}[section]
\newtheorem{lemma}[theorem]{Lemma}

\newtheorem{cy}[theorem]{Corollary}

\theoremstyle{definition}

\newtheorem{rk}[theorem]{Remark}

\newcounter{ppp}

\begin{document}

\title{On pairs of finitely generated subgroups in free groups.}
\author{A. Yu. Olshanskii \thanks{The
author was supported in part by the NSF grant DMS 1161294 and by the Russian Fund for Basic Research  grant 11-01-00945}}
\maketitle

\begin{abstract} We prove that for arbitrary two finitely generated subgroups
$A$ and $B$ having infinite index in a free group $F,$ there is a subgroup $H\le B$
with finite index $[B:H]$ such that the subgroup generated by $A$ and $H$ has infinite
index in $F$. The main corollary of this theorem says that a free group of  free rank 
$r\ge 2$ admits a faithful highly transitive action on an infinite set, whereas the restriction
of this action to any finitely generated subgroup of infinite index in $F$ has no infinite
orbits.
\end{abstract}

{\bf Key words:} free group, coset graph, highly transitive action

\medskip

{\bf AMS Mathematical Subject Classification: 20E05, 20B22, 20E07, 20E15, 54H15, 20F05}


\section{Introduction} 

The most known property of a pair of finitely generated subgroups $A$ and $B$
in a free group $F$ is the theorem of Howson (\cite{Ho}, \cite{LS}, I.3.13) saying
that the intersection $A\cap B$ is also finitely generated. The sharp estimate
of the free rank of $A\cap B$ in terms of (free) ranks of $A$ and $B$ was recently
published by Mineyev  \cite{M}. It confirms the old conjecture of Hanna Neumann  \cite{N}.

On the contrary, the subgroup $\langle A,B \rangle $ generated by $A\cup B$ has finite
rank for the obvious reason, and certainly this subgroup can have finite index in $F$ or can just coincide with $F$ when both $A$ and
$B$ are ''small'', that is, have infinite index in $F.$ However we prove in this paper
that $\langle A,H\rangle$ is still ``small'' for some subgroup $H$ which is
virtually equal to $B$. More precise formulation is given by

\begin{theorem}\label{2subgr} Let $A$ and $B$ be finitely generated subgroups of infinite index in 
a free group $F$. Then there is a subgroup $H\le B$ with finite index in $B$ such that
the subgroup $\langle A,H\rangle$ has infinite index in $F$.

In addition, for every finite subset $S\subset F\backslash A$, the subgroup $H$ can be chosen so that $\langle A,H\rangle\cap S =\emptyset.$
\end{theorem}

A particular case of Theorem \ref{2subgr} with a {\it cyclic} subgroup $B$ is contained in the theorem 5 of the paper \cite{BO}.  Theorem \ref{2subgr} is proved in Section \ref{main}; and in
Section \ref{next}, we formulate some corollaries obtained by the iterrated applications of this
theorem. The main corollary is

\begin{cy} \label{glav} Every noncyclic finitely generated free group $F$ (1) admits a 
highly transitive action $\circ$ on an infinite set $S,$ such that 
(2) the restriction of this action
to any finitely generated subgroup of infinite index in $F$ is locally finite.

Furthermore, every action of $F$ with properties (1) and (2) is faithful.
\end{cy}

Recall that an action of a group $G$ is {\it faithful} if its kernel is trivial. It is called {\it highly transitive} if it is $k$-transitive for every integer $k\ge 1.$ We call an action {\it locally finite} if all the orbits of this action are finite. 

Thus we have antipodal
properties (1) and (2) of the action $\circ$ of $F.$
On the one hand, an  easy observation shows that the restriction of an infinite transitive action with  the property(2) to a subgroup of {\it finite} index or to a nontrivial {\it normal} subgroup of $F$ cannot be locally finite (see Remark \ref{Rek}).
On the other hand the free group $F$ is saturated by finitely generated subgroups of infinite index. The $(r-1)$-generated
subgroups, where $r$ is the free rank of $F,$ constitute a minor part of this enormous family.

\section{Preliminaries: coset graphs, cores, and coverings} 

Let $F = F(X)$ be a free group with free basis $X$. If $H$ is a subgroup of $F$,
then the vertices of the {\it coset graph} $\Gamma=\Gamma(H,X)$
are the right cosets $Hg$ of the subgroup $H$ in $F$, and for every coset $Hg$ and
$x\in X^{\pm 1},$ we have an edge $e=e(Hg,x)$ with the initial vertex $e_-= Hg$ and the terminal vertex
$e_+=Hgx.$ The edge $e$ is labeled by $Lab(e)=x$, while the inverse edge $e^{-1}$
is labeled by $x^{-1},$ and $e^{-1}_-=Hgx$, $e^{-1}_+=Hg$.

Recall that if $e_i$ and $e_{i+1}$ are two consecutive edges in a path $p = e_1\dots e_n$ of (combinatorial) length $n$,
then $(e_{i+1})_-= (e_i)_+$ for $i=1,\dots,n-1$. The path $p$ of $\Gamma$  is {\it reduced} if $e_{i+1}\ne e_i^{-1}$ for $i=1,\dots,n-1$.
This is equivalent to saying that the word $Lab(p)=Lab(e_1)\dots Lab(e_n)$ is reduced. One defines
the initial and the terminal vertices of $p$ by the rules $p_-=(e_1)_-$ and $p_+=(e_n)_+.$ A path $p$  is {\it closed} if $p_-=p_+$. By the definition of inverse path, we have  $p^{-1}=e_n^{-1}\dots e_1^{-1}.$ For every vertex $v,$ there is a path $q$
of length $0$ with $q_-=q_+=v.$

Further, speaking on graphs labeled over some alphabet $X,$ we always
imply that for every edge $e$ we have a unique inverse $e^{-1}$ with $(e^{-1})^{-1}=e.$ We see  from the definition of the coset graph that it is connected, it has a base point
corresponding to the subgroup $H$, and every vertex $v$ of $\Gamma$ has a {\it standard } star $star(v)=star_{\Gamma}(v)$, i.e.,
the set of edges with the initial vertex $v$ has exactly one $x$-{\it edge} ($=$ the edge  labeled by $x$) for every $x\in X^{\pm 1}$. The following converse statement is known.

\begin{lemma}\label{cogr} Let $F=F(X)$ be a free group with free basis $X$. Then every connected labeled
graph $\Gamma$ with a basepoint, where each vertex $v$ is standard, is the coset graph of
a subgroup $H\le F$, and the base point corresponds to the coset $H.$
\end{lemma}

\proof Since every vertex $v$ is standard, for every $x\in X$, we have a 
bijection $f_x$ on the set of all vertices $V$: it maps each $v$ to the terminal vertex of the $x$-edge $e\in star(v)$. Since $F$ is free, the collection of all bijections $f_x$ extends
to an action $\circ$ of the entire $F$ on $V.$ 

The (right) action $\circ$ is transitive since $\Gamma$ is a connected graph. Let $H$ be the stabilizer of the base point $o\in\Gamma$. Since for $g,g'\in F$ and for any action $\circ$, one has  $o\circ g=o\circ g'$ if and only if $Hg=Hg'$, we have obtained the bijection $o\circ g\leftrightarrow Hg$, and the action $f_x$ of a free generator $x$ corresponds to the right multiplication by $x.$
\endproof

If $H\le H_1\le F$, then for $g,g'\in F$, the equality $Hg=Hg'$ implies $H_1g=H_1g'$, and therefore
the mapping $Hg\mapsto H_1g$ induces a label-preserving surjective graph morphism $\phi_H^{H_1}:\Gamma\to\Gamma_1$ of the corresponding coset graphs, preserving the base point. 
The restriction of the mapping $\phi_H^{H_1}$ to every star is bijective because all the labeled stars
are standard over the alphabet $X^{\pm 1}.$ In other words, $\phi_H^{H_1}$ is a {\it  covering} of
labeled graphs.
Obviously, the multiplicity of this covering (or its {\it index} $=$ the number of preimages
of any vertex or edge of $\Gamma_1$ in $\Gamma$) is equal to the number of
cosets of $H$ in any coset of $H_1$, i.e., it is equal to the index $[H_1:H]$ of the subgroup $H$ in $H_1$. We will apply the easy converse statement:

\begin{lemma}\label{cov} Let $\phi:\Gamma\to\Gamma_1$ be a preserving base point covering of index $j$ for
connected labeled graphs $\Gamma$ and $\Gamma_1$ with standard stars over the same alphabet $X^{\pm 1}$.
Then we have the inclusion $H\le H_1$ for the subgroups corresponded to these graphs in Lemma \ref{cogr}, and $[H_1:H]=j.$
\end{lemma}

\proof By Lemma \ref{cogr}, one may identify the graphs $\Gamma$ and $\Gamma_1$, respectively, with coset graphs of some subgroups $H$ and $H_1$ of the free group $F=F(X)$, where $o=H$ and $o_1=H_1$ are the base points. Since every star in $\Gamma$ is standard, for every word $w$ over the alphabet
$X^{\pm 1}$ and every vertex $v$ of $\Gamma$, there is a unique path $p$ with $p_-=v$ and $Lab(p)=w$.  It follows from the definition of coset graph that a path $p$ of $\Gamma$
originated at $o$ is closed if and only if $H\cdot Lab(p)=H$, i.e., iff  $Lab(p)$ represents  an element from the
subgroup $H.$ Since the projection $p_1$ of a closed path $p$ to $\Gamma_1$ is closed too, we have
$w=Lab(p)=Lab(p_1)\in H_1$ for every $w\in H$. Hence $H\le H_1$. Finally, the index of the covering is equal to $[H_1:H]$ as this was remarked before the lemma.
\endproof

Thus we see that a reduced word over the alphabet $X^{\pm 1}$ represents an element of $H$ if and only if $w$ is the
label of a reduced closed path of the coset graph $\Gamma$ with origin at the base point $o.$ So
one can consider the smallest subgraph $\cal C$ of $\Gamma$ containing $o$ and containing all the reduced closed paths
of $\Gamma$ originated at $o$. This labeled graph ${\cal C}=core(H,X)=core(H)$ with the base point $o$ is called the
{\it core} of $\Gamma.$

It follows from the definition that no vertex of the connected
graph $\cal C$, except for the base point $o,$ has degree $\le 1$ in $\cal C$. (Recall  that the {\it degree} of
a vertex $v$ in a graph $\cal C$ is the number of edges in $star_{\cal C}(v).$) If the degree of $o$ is $1$, then there is a 
unique maximal path $p$ in $\cal C$ such that $p_-=o$ and the degree of every its vertex in $\cal C$, 
except for the terminal vertex $p_+$, does not exceed $2$. Let us call this path the {\it handle} of
$\cal C.$  We will suppose that the handle of the core is of length $0$ if $\cal C$ has no
vertices of degree $1$. 

So the handle $p$ is attached at $p_+,$ to the part $\overline{\cal C}$ of
the core having no vertices of degree $1$. Removing the handle from $\cal C$ one
obtains the core $\overline{\cal C}=core(H')$ with base point $p_+$ for a conjugate subgroup $H'=Lab(p)^{-1}H (Lab(p))$. 

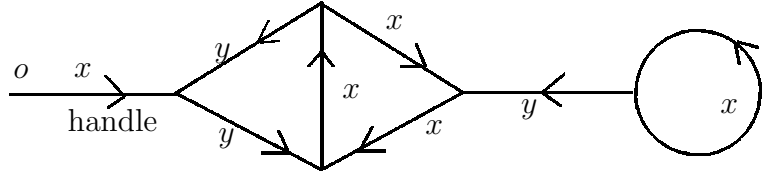
\begin{figure}[h!] \label{fig1}
\begin{center}
\unitlength 1mm 
\linethickness{0.4pt}
\ifx\plotpoint\undefined\newsavebox{\plotpoint}\fi 
\begin{picture}(104.026,39.374)(0,50)
\thicklines
\put(4.25,73.75){\line(1,0){21.75}}
\multiput(26,73.75)(.0542582418,.0336538462){364}{\line(1,0){.0542582418}}
\put(45.75,86){\line(0,-1){22.25}}
\multiput(45.75,63.75)(-.0625,.03125){8}{\line(-1,0){.0625}}
\multiput(45.25,64)(-.063973064,.0336700337){297}{\line(-1,0){.063973064}}
\put(26.25,74){\line(0,-1){.25}}
\put(26.25,73.75){\line(1,0){.25}}
\multiput(45.75,85.75)(.0537249284,-.0336676218){349}{\line(1,0){.0537249284}}
\put(64.5,74){\line(0,1){0}}
\multiput(64.5,74)(-.0616776316,-.0337171053){304}{\line(-1,0){.0616776316}}
\put(45.75,63.75){\line(0,1){.5}}
\put(64.25,74){\line(1,0){22.75}}
\put(104.026,74){\line(0,1){.4618}}
\put(104.014,74.462){\line(0,1){.4604}}
\put(103.975,74.922){\line(0,1){.4575}}
\put(103.911,75.38){\line(0,1){.4532}}
\put(103.821,75.833){\line(0,1){.4475}}
\multiput(103.706,76.28)(-.027924,.088073){5}{\line(0,1){.088073}}
\multiput(103.567,76.721)(-.032795,.086378){5}{\line(0,1){.086378}}
\multiput(103.403,77.153)(-.031303,.070345){6}{\line(0,1){.070345}}
\multiput(103.215,77.575)(-.030153,.058704){7}{\line(0,1){.058704}}
\multiput(103.004,77.986)(-.033382,.05693){7}{\line(0,1){.05693}}
\multiput(102.77,78.384)(-.031943,.048107){8}{\line(0,1){.048107}}
\multiput(102.514,78.769)(-.030735,.041111){9}{\line(0,1){.041111}}
\multiput(102.238,79.139)(-.032981,.039332){9}{\line(0,1){.039332}}
\multiput(101.941,79.493)(-.031612,.033687){10}{\line(0,1){.033687}}
\multiput(101.625,79.83)(-.033442,.031871){10}{\line(-1,0){.033442}}
\multiput(101.29,80.148)(-.039076,.033284){9}{\line(-1,0){.039076}}
\multiput(100.939,80.448)(-.040872,.031052){9}{\line(-1,0){.040872}}
\multiput(100.571,80.727)(-.047859,.032313){8}{\line(-1,0){.047859}}
\multiput(100.188,80.986)(-.049587,.029593){8}{\line(-1,0){.049587}}
\multiput(99.791,81.223)(-.05847,.030605){7}{\line(-1,0){.05847}}
\multiput(99.382,81.437)(-.070101,.031845){6}{\line(-1,0){.070101}}
\multiput(98.961,81.628)(-.086122,.033461){5}{\line(-1,0){.086122}}
\multiput(98.531,81.795)(-.087855,.028603){5}{\line(-1,0){.087855}}
\put(98.092,81.938){\line(-1,0){.4466}}
\put(97.645,82.057){\line(-1,0){.4525}}
\put(97.193,82.15){\line(-1,0){.457}}
\put(96.736,82.218){\line(-1,0){.46}}
\put(96.276,82.26){\line(-1,0){.9236}}
\put(95.352,82.267){\line(-1,0){.4606}}
\put(94.891,82.232){\line(-1,0){.458}}
\put(94.433,82.171){\line(-1,0){.4539}}
\put(93.98,82.085){\line(-1,0){.4483}}
\multiput(93.531,81.974)(-.088286,-.027244){5}{\line(-1,0){.088286}}
\multiput(93.09,81.837)(-.086628,-.032127){5}{\line(-1,0){.086628}}
\multiput(92.657,81.677)(-.070584,-.030759){6}{\line(-1,0){.070584}}
\multiput(92.233,81.492)(-.058935,-.029699){7}{\line(-1,0){.058935}}
\multiput(91.821,81.284)(-.057186,-.032941){7}{\line(-1,0){.057186}}
\multiput(91.42,81.054)(-.048352,-.031571){8}{\line(-1,0){.048352}}
\multiput(91.033,80.801)(-.041347,-.030417){9}{\line(-1,0){.041347}}
\multiput(90.661,80.527)(-.039585,-.032677){9}{\line(-1,0){.039585}}
\multiput(90.305,80.233)(-.03393,-.031351){10}{\line(-1,0){.03393}}
\multiput(89.966,79.92)(-.032128,-.033195){10}{\line(0,-1){.033195}}
\multiput(89.644,79.588)(-.033584,-.038818){9}{\line(0,-1){.038818}}
\multiput(89.342,79.238)(-.031366,-.040631){9}{\line(0,-1){.040631}}
\multiput(89.06,78.873)(-.032682,-.047608){8}{\line(0,-1){.047608}}
\multiput(88.798,78.492)(-.029974,-.049357){8}{\line(0,-1){.049357}}
\multiput(88.559,78.097)(-.031056,-.058232){7}{\line(0,-1){.058232}}
\multiput(88.341,77.689)(-.032385,-.069853){6}{\line(0,-1){.069853}}
\multiput(88.147,77.27)(-.028437,-.071551){6}{\line(0,-1){.071551}}
\multiput(87.976,76.841)(-.02928,-.087632){5}{\line(0,-1){.087632}}
\multiput(87.83,76.403)(-.03043,-.11141){4}{\line(0,-1){.11141}}
\put(87.708,75.957){\line(0,-1){.4517}}
\put(87.612,75.505){\line(0,-1){.4564}}
\put(87.54,75.049){\line(0,-1){.4597}}
\put(87.495,74.589){\line(0,-1){.9235}}
\put(87.48,73.666){\line(0,-1){.4609}}
\put(87.512,73.205){\line(0,-1){.4584}}
\put(87.569,72.747){\line(0,-1){.4545}}
\put(87.652,72.292){\line(0,-1){.4492}}
\multiput(87.76,71.843)(.0332,-.11062){4}{\line(0,-1){.11062}}
\multiput(87.892,71.4)(.031458,-.086874){5}{\line(0,-1){.086874}}
\multiput(88.05,70.966)(.030213,-.070819){6}{\line(0,-1){.070819}}
\multiput(88.231,70.541)(.029244,-.059163){7}{\line(0,-1){.059163}}
\multiput(88.436,70.127)(.032499,-.057439){7}{\line(0,-1){.057439}}
\multiput(88.663,69.725)(.031197,-.048594){8}{\line(0,-1){.048594}}
\multiput(88.913,69.336)(.030097,-.04158){9}{\line(0,-1){.04158}}
\multiput(89.184,68.962)(.03237,-.039836){9}{\line(0,-1){.039836}}
\multiput(89.475,68.603)(.031088,-.034171){10}{\line(0,-1){.034171}}
\multiput(89.786,68.262)(.032946,-.032383){10}{\line(1,0){.032946}}
\multiput(90.115,67.938)(.034702,-.030495){10}{\line(1,0){.034702}}
\multiput(90.462,67.633)(.040388,-.031679){9}{\line(1,0){.040388}}
\multiput(90.826,67.348)(.047354,-.033048){8}{\line(1,0){.047354}}
\multiput(91.205,67.083)(.049125,-.030354){8}{\line(1,0){.049125}}
\multiput(91.598,66.841)(.057991,-.031504){7}{\line(1,0){.057991}}
\multiput(92.004,66.62)(.069601,-.032923){6}{\line(1,0){.069601}}
\multiput(92.421,66.422)(.071329,-.028988){6}{\line(1,0){.071329}}
\multiput(92.849,66.249)(.087403,-.029956){5}{\line(1,0){.087403}}
\multiput(93.286,66.099)(.11117,-.03129){4}{\line(1,0){.11117}}
\put(93.731,65.974){\line(1,0){.451}}
\put(94.182,65.873){\line(1,0){.4559}}
\put(94.638,65.799){\line(1,0){.4593}}
\put(95.097,65.749){\line(1,0){.4614}}
\put(95.558,65.726){\line(1,0){.462}}
\put(96.02,65.728){\line(1,0){.4611}}
\put(96.481,65.756){\line(1,0){.4588}}
\put(96.94,65.81){\line(1,0){.4551}}
\put(97.395,65.889){\line(1,0){.45}}
\multiput(97.845,65.993)(.11087,.03235){4}{\line(1,0){.11087}}
\multiput(98.289,66.123)(.087114,.030786){5}{\line(1,0){.087114}}
\multiput(98.725,66.276)(.07105,.029666){6}{\line(1,0){.07105}}
\multiput(99.151,66.454)(.069284,.033584){6}{\line(1,0){.069284}}
\multiput(99.567,66.656)(.057688,.032055){7}{\line(1,0){.057688}}
\multiput(99.97,66.88)(.048833,.030821){8}{\line(1,0){.048833}}
\multiput(100.361,67.127)(.047038,.033497){8}{\line(1,0){.047038}}
\multiput(100.737,67.395)(.040085,.032062){9}{\line(1,0){.040085}}
\multiput(101.098,67.683)(.03441,.030824){10}{\line(1,0){.03441}}
\multiput(101.442,67.992)(.032637,.032696){10}{\line(0,1){.032696}}
\multiput(101.769,68.319)(.030762,.034466){10}{\line(0,1){.034466}}
\multiput(102.076,68.663)(.03199,.040143){9}{\line(0,1){.040143}}
\multiput(102.364,69.025)(.033413,.047098){8}{\line(0,1){.047098}}
\multiput(102.631,69.401)(.030733,.048889){8}{\line(0,1){.048889}}
\multiput(102.877,69.792)(.031951,.057746){7}{\line(0,1){.057746}}
\multiput(103.101,70.197)(.033459,.069345){6}{\line(0,1){.069345}}
\multiput(103.302,70.613)(.029538,.071104){6}{\line(0,1){.071104}}
\multiput(103.479,71.039)(.03063,.087169){5}{\line(0,1){.087169}}
\multiput(103.632,71.475)(.03215,.11093){4}{\line(0,1){.11093}}
\put(103.761,71.919){\line(0,1){.4502}}
\put(103.864,72.369){\line(0,1){.4553}}
\put(103.943,72.824){\line(0,1){.4589}}
\put(103.995,73.283){\line(0,1){.7166}}
\multiput(17.25,72)(.05,.03333333){45}{\line(1,0){.05}}
\multiput(19.5,73.5)(-.03358209,.03731343){67}{\line(0,1){.03731343}}
\multiput(17.25,76)(-.03125,.03125){8}{\line(0,1){.03125}}
\put(41.5,83.5){\line(0,1){0}}
\put(37.75,66){\line(1,0){3.75}}
\multiput(40,69)(.03289474,-.07236842){38}{\line(0,-1){.07236842}}
\multiput(44,77.5)(.03365385,.03846154){52}{\line(0,1){.03846154}}
\multiput(45.75,79.5)(.03333333,-.04444444){45}{\line(0,-1){.04444444}}
\put(56.25,77.5){\line(1,0){3.25}}
\multiput(59.5,77.5)(-.03289474,.07894737){38}{\line(0,1){.07894737}}
\put(50.5,66.25){\line(1,0){3.75}}
\multiput(50.75,66.25)(.03365385,.0625){52}{\line(0,1){.0625}}
\multiput(75,74)(.04,.03333333){75}{\line(1,0){.04}}
\multiput(75,74)(.04807692,-.03365385){52}{\line(1,0){.04807692}}
\multiput(100.75,80.5)(.0333333,-.1833333){15}{\line(0,-1){.1833333}}
\multiput(101,80.75)(.1195652,-.0326087){23}{\line(1,0){.1195652}}
\put(32,67.5){$y$}
\put(31.5,78.75){$y$}
\put(72.25,71.25){$y$}
\put(54.25,82.5){$x$}
\put(98.75,71.25){$x$}
\put(48.5,73.25){$x$}
\put(59.5,68.5){$x$}
\put(4.75,76){$o$}
\put(12.824,76.068){$x$}
\thinlines
\put(62.647,107.602){\line(1,0){.21}}
\put(26.278,119.164){\line(0,1){.21}}
\multiput(37.21,80.693)(.03344494,.06688989){44}{\line(0,1){.06688989}}
\multiput(37.63,80.693)(.360386,.030032){7}{\line(1,0){.360386}}
\put(11.983,68.92){handle}
\end{picture}

\end{center}
\caption{$core(H)$ for $H=\langle xyxyx^{-1}, xyx^3y^{-1}x^{-1},xyx^{-1}y^{-1}xyxy^{-1}x^{-1}\rangle\le F(x,y)$}
\end{figure}

Note that the notion of core of a graph was introduced 
by J.Stallings \cite{S}, but our definition is slightly different: Stallings' did not 
include the handle in the core and considered only the subgraph $\overline{\cal C},$ but we cannot
do this since we may not replace the subgroup by a conjugate one in the next section. 
Besides, we work only with labeled graphs.

Every reduced word representing an element
of $H$ can be read on a unique closed path starting at the base point. This path belongs to the minimal subgraph containing all the closed paths labeled by the reduced words generating $H.$
It follows that the core $\cal C$ is a finite graph if the subgroup $H$ is finitely generated.

A star $star_{\cal C}(v)$ can be smaller than the star $star_{\Gamma}(v)$ in a bigger graph. However 
if all the stars of $\Gamma$ are standard over some labeling alphabet $X$, then every subgraph $\cal E$
of $\Gamma$ is {\it regular}, i.e. every $star_{\cal E}(v)$
contains at most one $x$-edge for every $x\in X^{\pm 1}.$

Since all the edges of reduced closed paths of the coset graph $\Gamma$ with origin $o$ belong to  the core $\cal C$ and $\Gamma$ is a connected
graph, one obtains $\Gamma$ by attaching infinite labeled trees ${\cal T}_1$, ${\cal T}_2,\dots$ 
(``hairs'' in terminology of \cite{S}) to different vertices $v_1, v_2,\dots$ of $\cal C$ in such a way
that exactly one vertex of each tree ${\cal T}_i$ is identified with a vertex $v_i$ of $\cal C$. One attaches such a tree
to $v_i$ only if $v_i$ is a {\it deficit} vertex in the core, i.e. $star_{\cal C}(v_i)$
does not contain an $x$-edge for some $x\in X^{\pm 1}.$ 

In each of the attached trees, all stars $star_{{\cal T}_i}(v)$ of its vertices $v$ are standard, except for the root vertex $o({\cal T}_i),$ that
coincides with the  vertex $v_i$ of $\cal C$. Still the star $star(o({\cal T}_i))$ 
must be standard in the whole coset graph $\Gamma$; so the labeling of the edges of $star_{{\cal T}_i}(o({\cal T}_i))$ 
must compliment the labeling of the star $star_{\cal C}(v_i)$ of the deficit vertex $v_i$ in the core. Hence the regular graph 
$core(H)$ completely determines the coset graph $\Gamma$ of $H$ up to isomorphism. Furthermore, Lemma \ref{cogr} and  the above reconstruction of $\Gamma$ from the core by attaching of labeled trees
at deficit vertices, prove the following

\begin{lemma}\label{core}
Let $\cal C$ be a connected regular graph with a base point $o,$ labeled over an alphabet $X.$ If every vertex of $\cal C,$ except for $o,$ has degree at least $2$ in $\cal C,$ then $\cal C$ is equal to $core(H)$ with base point $o$ for a subgroup $H\le F(X).$
\end{lemma}
$\Box$

Now let $H$ be a finitely generated subgroup of $F.$ Since $core(H)$ is a finite graph, the coset 
graph $\Gamma=\Gamma(H)$ is finite if and only if it has no trees attached to the core. In other
words, the index of a finitely generated subgroup $H\le F$ is finite in $F$ if and only if $core(H)$ has no
deficit vertices.

\begin{lemma}\label{def} If a connected regular graph $\Delta$ labeled over an alphabet $X$ has a deficit vertex,
then every its nonempty subgraph $\cal E$ has its own deficit vertex. Moreover if the subgraph $\cal E$ is finite, 
 then either it has at least two distinct deficit vertices or it has a deficit vertex $v$ with deficit $>1$, i.e. there are two different letters $x,x'\in X^{\pm 1}$ such that $star_{\cal E}(v)$ has neither $x$- nor $x'$-edge. 
\end{lemma}
\proof If $\cal E$ is a proper subgraph with standard stars $star_{\cal E}(v)$, then there remain
no edges connecting $\cal E$ with the compliment $\Delta\backslash\cal E,$ contrary the connectedness
of $\Delta$. This proves the first assertion.

Assume that the subgraph $\cal E$ is finite and has only one deficit vertex $v$, and the deficit of $v$ in $\cal E$ is $1$. Since every standard star has even number of edges, we conclude that the sum of degrees of the vertices of $\cal E$ is odd. This contradicts the fact that the total number of edges in $\cal E$
is even since for every edge $e$ of $\cal E$, the edge $e^{-1}$ also belongs to $\cal E.$
The obtained contradiction proves the lemma.
\endproof

\section{Main lemmas and proof of Theorem \ref{2subgr}}\label{main}

Let $H$ be a subgroup 
of a free group $F=F(X)$ with a free base $X,$ and $\Gamma$ be the coset graph of $H$
with respect to $X$. 
For a subset $Y\subset X$, let $F(Y)$ denote the subgroup of $F$ generated by $Y.$
A vertex $v$ of a subgraph $E\subset \Gamma$ will be  called a $Y$-{\it deficit} vertex  in $E$ if there is $y\in Y^{\pm 1}$ such that $star_{E}(v)$ has no edge labeled by $y$.

We define the $Y$-frame $frame(H,Y)$ 
as the maximal connected subgraph of $core(H)$ containing the base point and having all edge labels in $Y^{\pm 1}$.

\begin{lemma}\label{one} If the subgroup $H$ is finitely generated,  has infinite index in $F,$ and it is not contained in $F(Y)$, then $H$ contains a subgroup $K$ 
of finite index in $H$ such that $core(K)$ has a deficit vertex which
does not belong to $frame(K,Y)$.
\end{lemma}

\proof  Assume first that the graph ${\cal C}=core(H)$ has an edge $e$ such that

(1) the removal of $e^{\pm 1}$ from $core(H)$ makes the remaining part $\cal C'$ of $\cal C$ disconnected
(let us call such an edge a {\it bridge} in $\cal C$) and

(2) $ Lab(e)\in X\backslash Y$.

Then the graph $\cal C'$ has two connected components $\cal D$ and $\cal E$, where
$\cal D$ contains the base point $o$ of $\cal C$, and so it contains the whole
$frame(H,Y)$ since $Lab(e)\notin Y^{\pm 1}$. 
By Lemma \ref{def}, the graph $\cal E$
must have a deficit vertex $v$ which remains to be deficit when one gets the edges $e^{\pm 1}$ back,
i.e., $v$ is deficit in $\cal C.$ 
It follows that in $core(H),$ the deficit vertex $v$ is separated from $frame(H,Y)$ by the edges $e^{\pm 1},$ and so one may choose $K=H$. 

Thus we may further assume that every bridge of $\cal C$
is labeled by a letter from $Y^{\pm 1}$.

Now we can find $K$ as a subgroup of $H$ with arbitrarily prescribed index $j\ge 2.$
To define $K$ we will construct a coverings $f:\Delta\to\Gamma$ of index $j$.

Let $\Gamma_0,\Gamma_1,\dots,\Gamma_{j-1}$ be copies of the coset graph 
$\Gamma=\Gamma(H)$ and ${\cal C}_i$ be the core of $\Gamma_i$ for $i=0,\dots,j-1$. 
Note that ${\cal C}_0$ has an edge $e_0$ labeled by a letter $x\in X\backslash Y$ because $H$ is not
a subgroup of $F(Y).$ The edge $e_0$ is not a bridge in ${\cal C}_0$ by the above conjecture.
Hence one obtains connected graphs after removal of $e_0$ from ${\cal C}_0$ or from $\Gamma_0$.

Let $e_0$  connect some vertices $u_0$ and $v_0$ of the graph ${\cal C}_0$.
By $e_i$, $u_i$, and $v_i$ we denote the copies of $e_0$, $u_0$, and $v_0,$ respectively,
in ${\cal C}_i$ ($i=1,\dots,j-1$). 

To define $\Delta$ we modify the disjoint union $\Gamma_0\sqcup\dots\sqcup \Gamma_{j-1}$ as
follows. Preserving all the vertices and all the edges except for $e_0,\dots,e_{j-1}$,
we redirect each $e_i.$ To be exact, we replace it by a new
edge $e'_i$ with the same label $x,$ but
$e'_i$ connects $u_i$ with $v_{i+1}$ (indices are taken modulo $j$). The 
edge $(e'_i)^{-1}$ is redirected respectively. 

The obtained graph $\Delta$ is connected
because for every $i$, the removal of the edges $e_i^{\pm 1}$ does not break the connectedness of $\Gamma_i$.

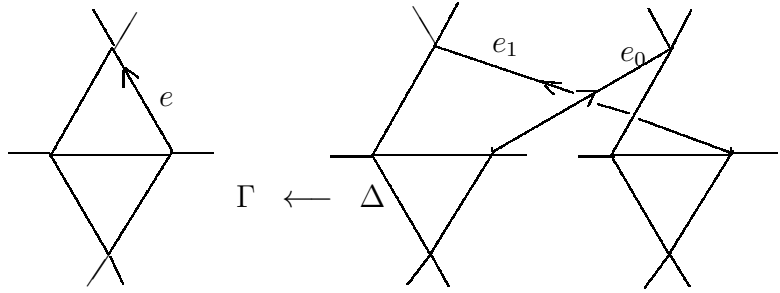
\begin{figure}[h!] \label{fig2}
\begin{center}
\unitlength 1mm 
\linethickness{0.4pt}
\ifx\plotpoint\undefined\newsavebox{\plotpoint}\fi 
\begin{picture}(1,63.5)(60,30)
\multiput(14.25,63)(.0336734694,.0591836735){245}{\line(0,1){.0591836735}}
\multiput(57,63)(.0336734694,.0591836735){245}{\line(0,1){.0591836735}}
\multiput(88.75,63)(.0336734694,.0591836735){245}{\line(0,1){.0591836735}}
\multiput(22.5,77.5)(.033613445,-.057773109){238}{\line(0,-1){.057773109}}
\multiput(30.5,63.75)(-.0337301587,-.0545634921){252}{\line(0,-1){.0545634921}}
\multiput(73.25,63.75)(-.0337301587,-.0545634921){252}{\line(0,-1){.0545634921}}
\multiput(105,63.75)(-.0337301587,-.0545634921){252}{\line(0,-1){.0545634921}}
\multiput(22,50)(-.033695652,.057608696){230}{\line(0,1){.057608696}}
\multiput(64.75,50)(-.033695652,.057608696){230}{\line(0,1){.057608696}}
\multiput(96.5,50)(-.033695652,.057608696){230}{\line(0,1){.057608696}}
\put(14.25,63.25){\line(0,1){0}}
\put(57,63.25){\line(0,1){0}}
\put(88.75,63.25){\line(0,1){0}}
\put(14.25,63.25){\line(1,0){16}}
\put(57,63.25){\line(1,0){16}}
\put(88.75,63.25){\line(1,0){16}}
\put(30.25,63.75){\line(0,1){.5}}
\put(73,63.75){\line(0,1){.5}}
\put(104.75,63.75){\line(0,1){.5}}
\put(24,74.5){\line(0,-1){2.25}}
\multiput(24.25,74.75)(.04605263,-.03289474){38}{\line(1,0){.04605263}}
\multiput(73.25,63.75)(.0582107843,.0337009804){408}{\line(1,0){.0582107843}}
\put(84.25,71.75){\line(1,0){2.5}}
\multiput(86.75,71.75)(-.0333333,-.0666667){30}{\line(0,-1){.0666667}}
\multiput(65.25,77.75)(.099462366,-.033602151){186}{\line(1,0){.099462366}}
\multiput(87,70)(.11184211,-.03289474){38}{\line(1,0){.11184211}}
\multiput(92.5,68)(.09141791,-.03358209){134}{\line(1,0){.09141791}}
\multiput(79.5,73)(.34375,-.03125){8}{\line(1,0){.34375}}
\put(30.25,63.5){\line(1,0){5.75}}
\put(8.75,63.5){\line(1,0){5.5}}
\put(22.75,77.25){\line(3,5){3}}
\multiput(22.75,77.75)(-.03353659,.05792683){82}{\line(0,1){.05792683}}
\put(22.25,50.25){\line(-2,-3){3}}
\multiput(22.25,49.75)(.03365385,-.07211538){52}{\line(0,-1){.07211538}}
\multiput(65.25,77.5)(.033505155,.06185567){97}{\line(0,1){.06185567}}
\put(65.5,78){\line(-3,5){3}}
\put(57.25,63){\line(-1,0){5.75}}
\put(72.75,63.25){\line(1,0){4.75}}
\multiput(65,50.25)(-.033505155,-.043814433){97}{\line(0,-1){.043814433}}
\multiput(64.75,50.25)(.03333333,-.06){75}{\line(0,-1){.06}}
\multiput(96.75,77.25)(.03353659,.06402439){82}{\line(0,1){.06402439}}
\multiput(96.75,77.25)(-.03370787,.06179775){89}{\line(0,1){.06179775}}
\put(84.5,63){\line(1,0){4.75}}
\put(104.5,63.5){\line(1,0){6.5}}
\multiput(96.75,50.25)(-.033653846,-.045673077){104}{\line(0,-1){.045673077}}
\multiput(96.5,50.25)(.03353659,-.05487805){82}{\line(0,-1){.05487805}}
\multiput(80,72.5)(.03289474,-.03289474){38}{\line(0,-1){.03289474}}
%
\put(90,75.75){$e_0$}
\put(28.75,69.75){$e$}
\put(73,77){$e_1$}
\put(39,56.25){$\Gamma\;\;\longleftarrow\;\;\Delta$}
\end{picture}

\end{center}
\caption{The covering $f: \Delta\to\Gamma$ for $j=2$}
\end{figure}

By definition, the function $f$ maps every vertex and every non-modified edge 
of each $\Gamma_i$ to its copy in $\Gamma$, and $f(e'_i)$ is the copy of $e_i$ in $\Gamma$.
Clearly, $f$ is a covering of $\Gamma$ of index $j$, and so by lemmas \ref{cogr} and \ref{cov}, $\Delta$ is the coset graph for a subgroup $K$ having index $j$ in $H.$ (In fact $K$ is normal in $H$
with cyclic quotients $H/K$.)

We choose the vertex $o_{\Delta}$ such that $f(o_{\Delta})=o_{\Gamma}$ and $o_{\Delta}$ belongs to ${\cal C}_0,$ as the base point of $\Delta$.
Let ${\cal C'}_i = {\cal C}_i\backslash\{e_i^{\pm 1}\}$ ($i=0,\dots, j-1$). Each graph ${\cal C'}_i$
is connected and contains the handle $p_i$ of ${\cal C}_i$ since $e_i$ is not a bridge in ${\cal C}_i$.
Hence the subgraph ${\cal C}_{\Delta}$ of
$\Delta$ formed by the edges $(e'_i)^{\pm 1}$ ($i=0,\dots,j-1$) and by the graphs ${\cal C'}_0$, ${\cal C'}_1,\dots,{\cal C'}_{j-1}$ 
without $j-1$ handles $p_1,\dots p_{j-1}$ is also connected. The graph ${\cal C}_{\Delta}$  has at most one vertex of degree $\le 1$ (must be equal to $(p_0)_-=o_{\Delta}$ if any exists), and ${\cal C}_{\Delta}$ contains all reduced closed
paths of $\Delta$ originated at $o_{\Delta}$ because every edge in $\Delta\backslash{\cal C}_{\Delta}$
belongs to some attached infinite tree. Hence ${\cal C}_{\Delta}$ is equal to $core(K)$.

The degrees of the vertices from ${\cal C'}_i$ in $core (K)$ do not exceed their degrees in $\Gamma_i$.
Now recall that $H$ is finitely generated and has infinite index in $F$; therefore $core(H)$ has a deficit vertex. So every ${\cal C}_i$ has a deficit vertex which remains to be deficit not only in ${\cal C'}_i$ but also in $\Delta$ by Lemma \ref{def}. So each of the parts ${\cal C'}_0,\dots,{\cal C'}_{j-1}$ of $\Delta$  has a deficit vertex in $\Delta.$

On the other hand, the subgraph $frame(K,Y)$ 
is contained in ${\cal C'}_0$ since ${\cal C'}_0$ is connected with the parts ${\cal C'}_1$ and ${\cal C'}_{j-1}$ by $x$-edges $(e'_0)^{\pm 1}$ and $(e'_{j-1})^{\pm 1}$ only, where $x\notin Y^{\pm 1}.$
Thus the deficit vertex of $\Delta$ belonging to ${\cal C'}_1$ is the desired one. 
\endproof

\begin{lemma}\label{two}
Let $A$ and $B$ be finitely generated subgroups of infinite
index in a free group $F$. Then there exist two subgroups $A_1$ and $B_1$ in $F$
such that

(a) $A_1$ is a subgroup of finite index in $A$;

(b) $B_1$ contains a subgroup $B_0$ having finite index in $B$;

(c) $B_1$ is finitely generated and contains $A_1$;

(d) the subgroup $B_1$ is of infinite index in $F$.
\end{lemma}

\proof Since $A$ is finitely generated, it is a free factor in a subgroup $E$
having finite index in $F$ by the theorem of M.Hall (see \cite{H}, \cite{LS}, I.3.10). Then $B\cap E$ has 
finite index in $B,$ and therefore it suffices to prove the lemma for the pair of subgroup
$(A, B\cap E)$ in the free group $E$ rather than for the pair $(A,B)$ in $F$. (Now and subsequently
we use that a subgroup of a free group is free and a subgroup of finite index in a finitely generated group is finitely
generated itself, see \cite{LS}, II,4.2.)  So one
may assume from the very beginning that $A$ is a free factor of $F$. In other
words, $F$ has a free basis $X$ such that $A$ is freely generated by a finite subset
$Y\subset X$. This subset is proper since  the subgroup $A$ is of infinite index in $F.$

If $B\le A=F(Y)$, then obviously one can put $A_1=B_1=A$ and $B_0=B$, which proves the lemma.
So we will assume that $B$ is not contained in $F(Y).$ By Lemma \ref{one}, $B$ has a
subgroup $K$ of finite index such that the subgraph $frame(K,Y))$ does not contain at least
one deficit vertex from $core (K)$. Again one may replace $B$ by
$K$ and prove the lemma constructing the pair $(A_1,B_1)$ for the pair $(A,K)$ rather than
for $(A,B)$ since $K$ has finite index in $B$. Therefore without change of our notation,
we may just assume further that $frame(B,Y)$ does not contain at least
one deficit vertex from $core(B).$

Now we are going to embed the graph $core(B)$ in the core of the coset graph of a bigger finitely generated
subgroup $B_1$. We will add edges to $core(B)$ as follows. Assume that $v$ is a $Y$-deficit vertex in $frame(B,Y)$. It is also $Y$-deficit in $core(B)$ by the definition of $Y$-frame. Then for some $y\in Y^{\pm 1}$, $core(B)$ has no edge terminating at $v$ and labeled by $y$. It follows that there exists 
a maximal path $q$ in the subgraph $frame(B,Y)$ (and in $core(B)$) such that $Lab(q)=y^n$ for some $n\ge 0$ and $q_-=v$.
Then we add a new edge going from $q_+$ to $q_-$ and labeled by $y.$
The maximality of $q$ implies that the extended graph is still regular. 

By Lemma \ref{core}, after each step of this procedure, we have a core for some subgroup of $F$. We will keep doing such extensions until the $Y$-frame has no $Y$-deficit vertices.  There appear
only finitely many new edges since we do not add new vertices. Therefore the 
subgroup $B_1$ given by the extended core is finitely generated, $core(B_1)$ contains $core(B)$,
and so $B_1\ge B$ . Since the subgraph $frame(B_1,Y)$ has no $Y$-deficit vertices, it is
equal to the core of the subgroup $A_1=F(Y)\cap B_1=A\cap B_1$, having 
finite index in $A$. (In fact, this argument also proves M.Hall's theorem:
$B\cap A$ is a free factor of a subgroup $A_1$ having finite index in $A.$)

Finally, $B_1$ has infinite index in $F$ since our extensions do not touch
the deficit vertex of $core(B)$ lying outside of $frame(B,Y)$. Hence
the properties (a) - (d) are obtained.
\endproof

For the next step, we need the following weaker form of Theorem \ref{2subgr}.

\begin{lemma}\label{cor} Let $A, B, C,\dots$ be a finite collection of finitely generated subgroups of infinite index in a free group $F$. Then there exist subgroups $A', B', C',\dots$ of finite
indices in $A, B, C,\dots$, respectively, such that the subgroup $\langle A',B',C',\dots \rangle$ has
infinite index in $F$.
\end{lemma}

\proof It follows from the properties (a) - (d) of Lemma \ref{two} that
one can put $A'=A_1$, $B'=B_0$ and obtain the subgroup $\langle A',B'\rangle\le B_1$ of infinite
index in $F$ with $A'$ and $B'$ having finite indices in $A$ and $B$, respectively. Therefore
the subgroups $A'$, $B'$, and $\langle A',B'\rangle$ are finitely generated.
In turn, there are a subgroup $H$ of finite index in
$\langle A',B'\rangle$ and a subgroup $C'$ of finite index in $C$  such that $\langle H,C'\rangle$ 
is of infinite index in $F$. But $H$ contains
some subgroups $A''$ and $B''$ of finite indices in $A'$ and $B'$, respectively,
and so in $A$ and in $B$. Therefore the subgroup $\langle A'',B'',C'\rangle$ also has
infinite index in $F$. Arguing in this way we complete the proof by induction on
the number of the subgroups in the finite set $\{A,B,C,\dots\}.$
\endproof

\begin{lemma}\label{three} Let $A$ and $B$ be finitely generated subgroups of infinite index in 
a free group $F$. Then there is a subgroup $B_2$ of $F$ such that

(a) $B_2$ is finitely generated;

(b) $B_2$ has infinite index in $F$;

(c) $B_2$ contains a subgroup $H$ of $B$ having finite index in $B$;

(d) the normalizer of $B_2$ in $F$ contains $A$.

\end{lemma}
\proof Let $A_1$ and $B_1$ be the subgroups of $F$ given by Lemma \ref{two}.
We choose some (finite) left transversal $\{a_1,\dots, a_m\}$ of $A_1$ in $A$
and define $L_i=a_iB_1a_i^{-1}$ ($i=1,...,m$). By Corollary \ref{cor}, there are
subgroups $N_1\le L_1,\dots, N_m\le L_m$ with finite indices $[L_i:N_i]$  such that 
the subgroup $B_2=\langle N_1,\dots N_m\rangle$ has infinite index in $F$. 

We have $N_i=a_iQ_ia_i^{-1}$, where $[B_1:Q_i]<\infty$. One can replace each $Q_i$
by a single $Q=\cap_{i=1}^m Q_i$, and $Q$ also has finite index in $B_1$. Moreover,
decreasing $Q$, we may assume that it is normal in $B_1$. Then it follows that $B_2$
is invariant under conjugations by any $a\in A$ since for any $a_i$ there are some $j$ 
and $a'\in A_1$ such that $aa_i=a_ja'$, and therefore
$$aN_ia^{-1}=aa_iQa_i^{-1}a^{-1}=a_ja'Qa'^{-1}a_j^{-1}=a_jQa_j^{-1}=N_j$$
Here $a'Qa'^{-1}=Q$ because $Q$ is normal in $B_1$ and $a'\in A_1\le B_1$ by Lemma \ref{two} (c).

Thus the properties (a), (b), and (d) of $B_2$ are proved. To complete the proof of the
lemma, we may assume that $a_1=1$, and so $B_2\ge N_1=Q$, but $Q$ is a subgroup of 
finite index in $B_1$, and therefore $Q$  must contain 
a subgroup $H$ of finite index in $B$ by Lemma \ref{two} (b). This proves the property (c).
\endproof

\medskip

{\bf Proof of Theorem \ref{2subgr}.} One may assume that $B\ne \{1\}.$ Let $B_2$ be a subgroup provided by Lemma \ref{three}. By Lemma \ref{three} (d), $B_2$ is normal in $\langle A,B_2\rangle=AB_2.$ It is finitely generated by Lemma \ref{three} (a). But a non-trivial finitely generated normal
subgroup of a free group must have finite index (see \cite{LS}, I.3.12), that is $|AB_2/B_2|<\infty.$ Hence $\langle A,B_2\rangle$ has infinite index in $F$ by Lemma \ref{three}(b). So has $\langle A,H\rangle$, where $H$ is a subgroup of $B\cap B_2$
provided by Lemma \ref{three} (c). The first claim of the theorem is proved.

To prove the second statement, we recall that given a finitely generated subgroup $A$ of a free group $F$
and a finite subset $S\subset F\backslash A,$ there is
a subgroup $M$ of finite index in $F$ such that $M\ge A$ and $M\cap S=\emptyset$ (see \cite{LS}, I.3.10).
Therefore we have $(M\cap\langle A,H\rangle)\cap S=\emptyset $ and consequently, 
$$\langle A,M\cap H\rangle\cap S=\langle M\cap A,M\cap H\rangle\cap S\le (M\cap\langle A,H\rangle)\cap S=\emptyset $$ So
to prove the second claim of the theorem we just replace $H$ by the subgroup $H'=H\cap M$
having finite indices in $H$ and in $B.$
$\Box$

\section{Corollaries for the actions of free groups}\label{next}

\begin{cy}\label{R} (a) A finitely generated free group $F$ has a subgroup $R$ of infinite index,
such that for every finitely generated subgroup $L\le F$ of infinite index in $F$, the
index $[L:L\cap R]$ is finite.

(b) Any subgroup $R$ with property (a) contains no non-trivial normal subgroups of $F.$
\end{cy}

\proof (a) Let $L_1,L_2,\dots$ be an enumeration of all finitely generated subgroups
having infinite index in $F$. Let us put $R_1=L_1$. For $i>1$, by induction, 
we define $R_i=\langle R_{i-1},H_i\rangle,$ where $H_i$ is a subgroup of finite index in $L_i$ 
such thar $R_i$ is a finitely generated subgroup of infinite index in $F$. Such a choice 
of $H_i$ is possible by Theorem \ref{2subgr} since by the inductive hypothesis, 
$R_{i-1}$ is a finitely generated subgroup of infinite index in $F$.

Now we define $R=\cup_{i=1}^{\infty}R_i$. We see that $R$ is a union of the members 
of the increasing series of the subgroups $R_i$-s having infinite index in $F.$ Hence 
the subgroup $R$ has infinite index in $F$ itself 
because $F$ is finitely generated, and so is every subgroup of finite index in $F$.
Since for every $i$, the subgroup $R$ contains a subgroup $H_i$ of finite index in $L_i,$
 the statement (a) is proved.
 
 (b) To prove the second statement we need one more lemma.
 \begin{lemma}\label{H} For every nontrivial normal subgroup $N$ of a free group $F=F(x_1,\dots,x_r)$ of finite rank $r\ge 2$, there exists an $r$-generated subgroup $H$ of infinite index in $F$ such that $F=HN$.
 \end{lemma}
 \proof We will apply a small cancellation argument. 
 Let $w$ be a non-trivial reduced word representing an element of $N$. Its length is
 denoted by $|w|$.
 There exist $2r$ positive words $u_1,\dots, u_{2r}$ in the alphabet $X$
 such that each of them has length at least $10|w|,$ and arbitrary subword of length $\ge |u_i|/10$
 of any $u_i$ does not occur as a subword in $u_j$ 
 for $j\ne i$ and has a unique occurrence  in $u_i$ (see the examples in \cite{LS}, V.10).
 
 Let $v_i$ be the reduced form of the product $u_{2i-1}wu_{2i-1}^{-1}x_iu_{2i}wu_{2i}^{-1}$
 ($i=1,\dots,r$).
 Clearly, $v_i \equiv x_i$ modulo $N$, and so $F=NH$ for the subgroup $H=\langle v_1,\dots,v_r\rangle.$
 It easily follows from the choice of the words $u_j$-s that the reduced form of the
 product $u_{2i-1}^{-1}x_iu_{2i}$ starts (respectively, ends) with a prefix of $u_{2i-1}^{-1}$ (with
 a suffix of $u_{2i}$) of length at least $\frac45|u_{2i-1}|$ (at least $\frac45|u_{2i}|).$
 Similarly, $v_i$ starts (ends) with a prefix of $u_{2i-1}$ (with
 a suffix of $u_{2i}^{-1}$) of length at least $\frac45|u_{2i-1}|$ (at least $\frac45|u_{2i}|).$
 Therefore if two words $v_i^{\pm 1}$ and $v_j^{\pm 1}$ are not mutually inverse, then
 the cancellations in their products affect less than $1/4$ of the letters of each of them.
 Hence the subgroup $H$ is freely generated by the words $v_1,...,v_r.$
 
 The same small cancellation argument shows that no product of the generators of $H$ is
 equal to one of the free generators $x_1,...,x_r$ of $F$, that is $H$ is a proper subgroup of $H.$
 Since by Schreier's formula (see \cite{LS},I.3.9), the free rank of every proper subgroup of 
 finite index in $F$ is strictly greater than $r,$ the subgroup $H$ must have infinite index
 in $F$, and the lemma is proved.
 \endproof 

To finish the proof of the statement (b) we may assume that the free rank $r$ of $F$ is at least $2$.
Assume that the subgroup $R$ contains a nontrivial subgroup $N$ normal in $F$. Let $H$ be the finitely
generated subgroup of infinite index in $F$ provided by Lemma \ref{H}. Since $R$ contains 
$N$ and, by (a), it contains a subgroup of finite index from $H$, we conclude that $R$ has finite index in $F$ because
$F=HN.$ The obtained contradiction with the property (a) completes the proof of (b). 
\endproof

\begin{rk}\label{Rek}  Obviously, the index $[L:L\cap R]$ is infinite for every subgroup $L$
of finite index. Also one cannot omit the assumption that $L$ is finitely generated
in part (a) of Corollary \ref{R}; in particular, it follows from (b) that the index 
$[L:L\cap R]$ is infinite for every nontrivial normal subgroup $L$ of $F$.
Indeed, otherwise there is a characteristic subgroup $N$ of $L$ such that $N\le L\cap R\le R$
and the factor group $L/N$ has finite exponent equal to $[L:L\cap R]$. Therefore $N$ is
a nontrivial normal subgroup of $F$ contrary to Lemma \ref{R} (b).
\end{rk}

The following reformulation of Corollary \ref{R} in terms of the actions of $F$
answers the question raised in
Open Problem 3 of \cite{BO}.

\begin{cy}\label{lf} (a) Every nontrivial finitely generated free group $F$ admits a 
transitive action on an infinite set $S$ such that the restriction of this action to
any finitely generated subgroup $L$ of infinite index in $F$ is locally finite.

(b) Every action $\circ$ of $F$ satisfying the conditions from (a) is faithful.
\end{cy}

\proof (a) Let $F$ act by right translations on the set $S$ of right cosets of the subgroup $R$
provided by Corollary \ref{R}. Then $S$ is infinite and the action on $F$ is transitive.
An $L$-orbit $(Rg)L$ of a point $Rg$ with respect to the action of a finitely generated subgroup $L$ with infinite $[F:L]$ has the same size as the orbit $RL'$ of the action of the conjugate subgroup $L'=gLg^{-1}$, which is also finitely generated and has infinite index in $F$. But the stabilizer
of the point $R$ under the action of $L'$ is $L'\cap R$. So the size of the orbit $RL'$ is equal
to the index $[L':L'\cap R]$, which is  finite by Corollary \ref{R}. The local finiteness is proved.

(b) The stabilizer $R$ of any point $v$ must have infinite index in $F$ since $F$ transitively 
acts on an infinite set. Since the restriction of the action $\circ$ to every finitely generated subgroup
$H$ of infinite index is locally finite, the $H$-orbit $v\circ H$ is finite. Since the size of this
orbit is equal to $[H:H\cap R]$, we conclude by Corollary \ref{R}, that $R$ does not contain nontrivial
normal subgroups of $F$. Hence the kernel of the action $\circ$ is trivial, that is this action is
faithful.
\endproof 

The approach used in the proof of corollaries \ref{R}, \ref{lf} can be employed to unite
different extremal properties of group actions. We combine two opposite properties
of actions, namely, locally finiteness and multiple transitivity.
First faithful highly 
transitive actions of non-cyclic free groups of finite ranks on infinite sets were
obtained in \cite{M} (see other proofs in \cite{D}, \cite{BO}). Similar results were proved for
free products of groups in \cite{GM}, \cite{Hi}, and \cite{G} (see also \cite{MS}). Quite recently
this theorem was extended to surface groups \cite{K}, to $Out(F_n)$ \cite{GG}, and finally,
to all non-elementary hyperbolic groups with trivial finite radical \cite{C}. Applying Theorem \ref{2subgr} and a theorem from \cite{BO}, one can sharpen Corollary \ref{lf} and obtain Corollary \ref{glav}.
\medskip

{\bf Proof of Corollary \ref{glav}.} As in the proof of Corollary \ref{lf}, the group $F$ acts by right multiplications
on the set of the right cosets of a subgroup $R.$ However the construction of $R$ is modified
in comparison with Corollary \ref{R} as follows. We will enumerate both finitely generated
subgroups of infinite index in $F$ and all $2k$-tuples $(g_1,\dots,g_k,g'_1,...,g'_k)$ of 
elements from $F$ for all $k\ge 1$. In the inductive definition of the increasing series $R_1\le R_2\le\dots$,
we introduce $R_i$ using the rule from Corollary \ref{R} if $i$ is odd. For even $i$-s, we
will apply the following theorem 6 from \cite{BO}:

\medskip

{\it Let $H$ be a finitely generated subgroup of infinite index in a free
group $F$ of rank $r > 1$. Let $(g_1H,\dots, g_kH)$ and $(g'_1H,\dots, g'_kH)$
($k\ge 1$) be
two $k$-tuples of pairwise different cosets. Then there
is a finitely generated subgroup $H'$ of infinite index in $F$ 
and an element $b\in F$, such that $H\le H'$ and $H'g_jb = H'g'_j$
for every $j=1,\dots,k$.}

\medskip

Now if $i$ is even, we take the first tuple $(g_1,\dots,g_k,g'_1,...,g'_k)$ (if any exists)
 in our enumeration such that (1) it  was not considered at the previous steps, (2)
the cosets $R_{i-1}g_1,\dots, R_{i-1}g_k$ are pairwise different, and (3)
the cosets $R_{i-1}g'_1,\dots, R_{i-1}g'_k$ are pairwise different too.
Then we apply the cited theorem from \cite{BO} to the subgroup $H=R_{i-1}$ and to the
elements $g_1,\dots g_k,g'_1,\dots, g'_k$ and obtain $R_i=H'$.
If there are no tuples with the properties (1)--(3) then we set $R_i=R_{i-1}.$

Again the subgroup $R=\cup_{i=1}^{\infty}R_i$ has infinite index and the properties
from Corollary \ref{lf} hold. In addition, the action is now $k$-transitive for any $k\ge 1$ because
if we have two $k$-tuples of pairwise distinct cosets $(Rg_1,...,Rg_k)$ and 
$(Rg'_1,\dots,Rg'_k)$, then the corresponding cosets remain pairwise distinct modulo every subgroup $R_{i-1}\le R$,
and so by the above construction, we should have $R_ig_jb = R_ig'_j$ for some integer $i$, $b\in F$, and $j=1,\dots,k$.
It follows that $Rg_jb = Rg'_j$, as required for the $k$-transitivity.
$\Box$

\medskip

\begin{rk} 
One can derive an effective bound for the numbers of edges and vertices in the coset
graph of the subgroup $H$ in terms of the coset graphs of the subgroups $A$ and $B$ 
given in the formulation of Theorem \ref{2subgr}; and for given $A$ and $B$, the set of generators of $H$ can be found
algorithmically. Taking into account the decidability of the membership problem for
finitely generated subgroups of free groups, one can construct a Turing machine enumerating
the generators of the subgroup $R$ in  Theorem \ref{2subgr} and its corollaries. Furthermore,
the last sentence of Theorem \ref{2subgr} is a clue to the machine enumeration of the compliment $F\backslash R.$
So the subgroup $R$ can be recursive in the statements of this paper.

However the
upper bounds for the parameters of $H$ based on our proof of  Theorem \ref{2subgr} are too rough and unsatisfactory.
The problem of finding of realistic estimates is open.

Another open problem was formulated in \cite{BO}: Can the action with the property (a)
from Corollary \ref{lf} have maximal growth ? (See details in \cite{BO}.)
\end{rk}

\medskip

{\bf Acknowledgements.} The author is thankful to Yuri Bahturin and Jonny Lomond for discussions.

\bigskip

{\bf Alexander A. Olshanskii:} Department of Mathematics, Vanderbilt University, Nashville
37240, U.S.A.

E-mail: alexander.olshanskiy@vanderbilt.edu


\begin{thebibliography}{99}
\label{bibl}
\newcommand{\bi}{\bibitem}

\bi{BO} Bahturin Yu. A. and Olshanskii A.Yu., Actions of maximal growth, 
Proceedings  of the London Math. Soc. 101 (2010), no.1, 27--72.

\bi{C} Chaynikov V.V., Properties of hyperbolic groups: free normal subgroups,
quasiconvex subgroups, and actions of maximal growth, Ph. D. Thesis, Vanderbilt
University (2012), available at http://etd.library.vanderbilt.edu/available/etd-06212012-172048/unrestricted/CHAYNIKOV.pdf .

\bi{D} Dixon J.D., Most finitely generated permutation groups are free, Bull.
London Math. Soc. 22 (1990), no.3, 222--226.

\bi{GG} Garion S. and Glasner Y., Highly transitive actions of $Out(F_n)$,
Groups, Geometry and Dynamics, 7 (2013), no. 2, 357--376.

\bi{GM} Glass A.M. and McCleary S.H., Highly transitive representations of free groups
and free products, Bull.Austral.Math.Soc. 43(1991), 19--36.

\bi{G} Gunhouse S.V., Highly transitive representations of free products on the natural
numbers, Arch. Math. 58 (1992), 435--443.

\bi{H} Hall, Marshall Jr., Coset representations in free groups, Trans.
Amer. Math. Soc. 67 (1949), 421 -- 432.

\bi{Hi} Hichin K.K., Highly transitive Jordan representations of free products.
J. London Math. Soc., 46 (1992), 81--91.

\bi{Ho} Howson A.G., On the intersection of finitely generated free groups,
J. London Math. Soc., 29 (1954), 428-434.


\bi{K} Kitroser D., Highly transitive actions of surface groups, Proc. Amer.
Math. Soc., 140 (2012), 3365--3375.

\bi{LS} Lyndon, R. and Schupp, P., Combinatorial Group Theory, Springer-Verlag,
2001.

\bi{M} McDonough T.P., A permutation representation of a free group, Quart. J.
Math. Oxford Ser. (2), 22 (1977), no.111, 353-356.

\bi{M} Mineyev I., Submultiplicativity and Hanna Neumann Conjecture, Ann. of Math.,
175 (2012), no. 1, 393-414.

\bi{MS} Moon S. and Stalder Y., Highly transitive actions of free products,
Algebr. Geom. Topol., 13 (2013), no. 1, 589--607.

\bi{N} Neumann H., On the intersection of finitely generated free groups. Publ. Math.
Debrecen, 4 (1956), 36-39; addendum, 5 (1957), 128.

\bi{S} Stallings, J.R., Topology of finite graphs, Invent. Math., 71 (1983), 551 - 565.

\end{thebibliography}
\end{document}